\documentclass[12pt]{article}
\usepackage{graphicx}
\font\bg=cmbx10 scaled\magstep1
\font\Bg=cmbx12 scaled\magstep3
\font\small=cmr8

\newtheorem{newlemma}{{\bf Lemma}}

\newenvironment{lemma}{\begin{newlemma}{\hspace{-0.5
em}{\bf.}}}{\end{newlemma}}

\newtheorem{newteorem}{{\bf Theorem}}

\newenvironment{theorem}{\begin{newteorem}{\hspace{-0.5
em}{\bf.}}}{\end{newteorem}}

\newtheorem{newkorolari}{{\bf Corollary}}

\newtheorem{newdefine}{{\bf Definition}}

\newenvironment{define}{\begin{newdefine}{\hspace{-0.5
em}{\bf.}}}{\end{newdefine}}

\newtheorem{newquestion}{{\bf Question}}

\newtheorem{newkonjek}{{\bf Conjecture}}

\newtheorem{newexample}{{\bf Example}}

\begin{document}
\tolerance=10000
\baselineskip18truept
\newbox\thebox
\global\setbox\thebox=\vbox to 0.2truecm{\hsize
0.15truecm\noindent\hfill}
\def\boxit#1{\vbox{\hrule\hbox{\vrule\kern0pt
     \vbox{\kern0pt#1\kern0pt}\kern0pt\vrule}\hrule}}
\def\qed{\lower0.1cm\hbox{\noindent \boxit{\copy\thebox}}\bigskip}
\def\ss{\smallskip}
\def\ms{\medskip}
\def\bs{\bigskip}
\def\c{\centerline}
\def\nt{\noindent}
\def\ul{\underline}
\def\ol{\overline}
\def\lc{\lceil}
\def\rc{\rceil}
\def\lf{\lfloor}
\def\rf{\rfloor}
\def\ov{\over}
\def\t{\tau}
\def\th{\theta}
\def\k{\kappa}
\def\l{\lambda}
\def\L{\Lambda}
\def\g{\gamma}
\def\d{\delta}
\def\D{\Delta}
\def\e{\epsilon}
\def\lg{\langle}
\def\rg{\rangle}
\def\p{\prime}
\def\sg{\sigma}
\def\ch{\choose}

\newcommand{\ben}{\begin{enumerate}}
\newcommand{\een}{\end{enumerate}}
\newcommand{\bit}{\begin{itemize}}
\newcommand{\eit}{\end{itemize}}
\newcommand{\bea}{\begin{eqnarray*}}
\newcommand{\eea}{\end{eqnarray*}}
\newcommand{\bear}{\begin{eqnarray}}
\newcommand{\eear}{\end{eqnarray}}
\centerline{\Bg}

\centerline {\Large{\bf Dominating Sets and Domination Polynomials
of Cycles}}
\bigskip
\baselineskip12truept \centerline{\bg Saeid
Alikhani$^{a,b,}${}\footnote{\baselineskip12truept\it\small
Corresponding author. E-mail: alikhani@yazduni.ac.ir} and Yee-hock
Peng$^{b,c}$} \baselineskip20truept \centerline{\it
$^{a}$Department of Mathematics} \vskip-8truept \centerline{\it
Yazd University} \vskip-8truept \centerline{\it 89195-741, Yazd,
Iran} \vskip-5truept \centerline{\it $^{b}$Institute for
Mathematical Research, and} \vskip-9truept \centerline{\it
$^{c}$Department of Mathematics,} \vskip-8truept \centerline{\it
University Putra Malaysia, 43400 UPM Serdang, Malaysia}
\vskip-8truept \centerline{\it } \vskip-0.2truecm
\nt\rule{16cm}{0.1mm}

\nt{\bg ABSTRACT}
\medskip

\nt{\it Let $G=(V,E)$ be a simple graph. A set $S\subseteq V$ is a
dominating set of $G$, if  every vertex in $V\backslash S$ is
adjacent to at least one vertex in $S$. Let ${\cal C}_n^i$ be the
family of dominating sets of a cycle $C_n$ with cardinality $i$,
and let $d(C_n,i)=|{\cal C}_n^i|$. In this paper,  we construct
${\cal C}_n^i$, and obtain a recursive formula for $d(C_n,i)$.
Using this recursive formula, we consider the polynomial
$D(C_n,x)=\sum_{i=\lceil\frac{n}{3}\rceil}^{n} d(C_n,i) x^{i}$,
which we call domination polynomial of cycles and obtain some
properties of this polynomial. }

\nt{ Keywords:} {\small  Dominating sets; Domination Polynomial;
Recursive formula; Cycle}


\nt\rule{16cm}{0.1mm}

\section{Introduction}
Let $G=(V,E)$ be a simple graph of order $|V|=n$.
 A set $S\subseteq V$ is a {\it dominating set} of $G$, if  every vertex in $V\backslash S$ is adjacent to at least one vertex in $S$.
The {\it domination number} $\gamma(G)$ is the minimum cardinality of a dominating set in $G$.
 For a detailed treatment of this parameter, the reader is referred to~\cite{domination}.
It is well known and generally accepted that the problem of
determining the dominating sets of an arbitrary graph is a
difficult one (see ~\cite{johnson}). Let ${\cal C}_n^i$ be the
family of dominating sets of a cycle $C_n$ with cardinality $i$
and let $d(C_n,i)=|{\cal C}_n^i|$. We call the polynomial
$D(C_n,x)=\sum_{i=\lceil\frac{n}{3}\rceil}^{n} d(C_n,i) x^{i}$,
the {\it domination polynomial of cycle}. For a detailed treatment
of domination polynomial of a graph, the reader is  referred
to~\cite{saeid}.

\ms

\nt    In the next section  we construct  the families of dominating sets of $C_n$ with cardinality $i$ by the families of dominating sets of $C_{n-1}, C_{n-2}$ and $C_{n-3}$ with cardinality $i-1$.
    We investigate the domination polynomial of cycle in Section 3.

\nt As usual we use $\lceil x \rceil$,  for the smallest integer greater than or equal to $x$.
 In this paper  we denote the set $\{1,2,...,n\}$ simply by $[n]$.

\section{Dominating sets of cycles}
 \nt       Let $C_n, n\geq 3$, be the cycle with $n$ vertices $V(C_n)=[n]$ and $E(C_n)=\{(1,2),(2,3),...,(n-1,n),(n,1)\}$.  Let $\mathcal{C}_n^i$ be the family of dominating sets of $C_n$ with cardinality $i$. We shall investigate dominating sets of cycles.
  A {\it simple path} is a path where all its internal vertices have degree two. We need the following lemmas to prove our main results in this section:

\begin{lemma}\label{lemma1}
The following properties hold for cycles,

\begin{enumerate}

\item[(i)]
(~\cite{gray},p.364)
$\gamma(C_{n})=\lceil\frac{n}{3}\rceil$.

\item[(ii)]
$\mathcal{C}_{j}^{i}=\emptyset$, if and only if $i>j$ or $i<\lceil\frac{j}{3}\rceil$. (by $(i)$ above).

\item[(iii)]
If a graph $G$ contains a simple path of length $3k-1$, then every dominating set of $G$ must contain at least $k$ vertices of the path. (by  observation).\quad\qed

\end{enumerate}
\end{lemma}

\nt To find a dominating set of $C_n$ with cardinality $i$, we do not need to consider dominating
sets of $C_{n-4}$  and $C_{n-5}$ with cardinality $i-1$. We show this in Lemma~\ref{lemma2}.
Therefore, we
only need to consider $\mathcal{C}_{n-1}^{i-1}, \mathcal{C}_{n-2}^{i-1}$, and $\mathcal{C}_{n-3}^{i-1}$.
The families of these dominating sets can be
empty or otherwise. Thus, we have eight combinations of whether these three families
are empty or not. Two of these combinations are not possible (see Lemma~\ref{lemma3}$(i)$ and $(ii)$).
Also, the combination that $\mathcal{C}_{n-1}^{i-1}=\mathcal{C}_{n-2}^{i-1}=\mathcal{C}_{n-3}^{i-1}=\emptyset$;
 no need to be considered because it
implies $\mathcal{C}_n^i=\emptyset$ (see Lemma~\ref{lemma3}$(iii)$). Thus we only need to consider five combinations or
cases. We consider this in Theorem~\ref{theorem1}.

\vspace{.2cm}

\begin{lemma}\label{lemma2}
If $Y$ is in $\mathcal{C}_{n-4}^{i-1}$ or $\mathcal{C}_{n-5}^{i-1}$  such that $Y\cup \{x\}\in \mathcal{C}_n^i$ for some $x\in [n]$, then $Y\in \mathcal{C}_{n-3}^{i-1}$.
\end{lemma}
\vspace{.1cm}
{\bf Proof.} Let $Y\in \mathcal{C}_{n-4}^{i-1}$ and $Y\cup \{x\}\in \mathcal{C}_n^i$ for some $x\in [n]$.
This means, by Lemma~\ref{lemma3}, we only need to consider $\{1,n-4\},\{2,n-4\}$ and $\{1,n-5\}$ as a subset of $Y$. In each  case, $Y\in \mathcal{C}_{n-3}^{i-1}$.
Now suppose that $Y\in \mathcal{C}_{n-5}^{i-1}$ and $Y\cup \{x\}\in \mathcal{C}_n^i$ for some $x\in [n]$.
This means, by Lemma~\ref{lemma3}, $\{1,n-5\}$ must be a subset of $Y$. So $Y\in \mathcal{C}_{n-3}^{i-1}$. \quad\qed

\vspace{.15cm}

\nt The following lemma follows from Lemma~\ref{lemma1}$(ii)$.

\vspace{.2cm}

\begin{lemma}\label{lemma3}
\begin{enumerate}
\item[(i)]
 If $\mathcal{C}_{n-1}^{i-1}=\mathcal{C}_{n-3}^{i-1}=\emptyset$, then $\mathcal{C}_{n-2}^{i-1}=\emptyset$,
\item[(ii)]
 If $\mathcal{C}_{n-1}^{i-1}\neq\emptyset$ and $\mathcal{C}_{n-3}^{i-1}\neq\emptyset$, then $\mathcal{C}_{n-2}^{i-1}\neq\emptyset$,
\item[(iii)]
 If $\mathcal{C}_{n-1}^{i-1}=\mathcal{C}_{n-2}^{i-1}=\mathcal{C}_{n-3}^{i-1}=\emptyset$, then $\mathcal{C}_{n}^{i}=\emptyset$.\quad\qed
\end{enumerate}
\end{lemma}

\vspace{.19cm}
\nt The following lemma follow from Lemma~\ref{lemma1}$(ii)$.

\vspace{.2cm}

\begin{lemma}\label{lemma4} Suppose that $\mathcal{C}_n^i\neq\emptyset$, then we have
\begin{enumerate}
\item[(i)]
 $\mathcal{C}_{n-1}^{i-1}=\mathcal{C}_{n-2}^{i-1}=\emptyset$, and $\mathcal{C}_{n-3}^{i-1}\neq \emptyset$
 if and only if $n=3k$ and $i=k$ for some $k\in N,$
 \item[(ii)]
  $\mathcal{C}_{n-2}^{i-1}=\mathcal{C}_{n-3}^{i-1}=\emptyset$ and $\mathcal{C}_{n-1}^{i-1}\neq \emptyset$
  if and only if $i=n$,
  \item[(iii)]
   $\mathcal{C}_{n-1}^{i-1}=\emptyset,\mathcal{C}_{n-2}^{i-1}\neq\emptyset$ and $ \mathcal{C}_{n-3}^{i-1}\neq \emptyset$
    if and only if $n=3k+2$ and $i=\lceil\frac{3k+2}{3}\rceil$ for some $k\in N$,
  \item[(iv)]
   $\mathcal{C}_{n-1}^{i-1}\neq\emptyset, \mathcal{C}_{n-2}^{i-1}\neq \emptyset$ and $\mathcal{C}_{n-3}^{i-1}= \emptyset$
   if and only if $i=n-1$,
  \item[(v)]
   $\mathcal{C}_{n-1}^{i-1}\neq\emptyset, \mathcal{C}_{n-2}^{i-1}\neq\emptyset$ and $\mathcal{C}_{n-3}^{i-1}\neq \emptyset$
    if and only if $\lceil\frac{n-1}{3}\rceil+1 \leq i \leq n-2$.
  \end{enumerate}
\end{lemma}
\nt{\bf Proof.}
\begin{enumerate}
\item[(i)]
 ($\Rightarrow$) Since ${\cal C}_{n-1}^{i-1}={\cal C}_{n-2}^{i-1}=\emptyset$, by Lemma~\ref{lemma1}$(ii)$,
 we have $i-1>n-1$ or $i-1<\lceil\frac{n-2}{3}\rceil$. If $i-1>n-1$, then $i>n$, and by Lemma~\ref{lemma1}$(ii)$, ${\cal C}_{n}^i=\emptyset$, a contradiction.
 So we have $ i< \lceil\frac{n-2}{3}\rceil +1 $, and since ${\cal C}_n^i\neq \emptyset$, together we have $\lceil\frac{n}{3}\rceil \leq i <\lceil\frac{n-2}{3}\rceil+1$, which give us $n=3k$ and $i=k$ for some $k\in N$.\\
($\Leftarrow$) If $n=3k$ and $i=k$ for some $k\in N$, then by Lemma~\ref{lemma1}$(ii)$,
we have ${\cal C}_{n-1}^{i-1}={\cal C}_{n-2}^{i-1}=\emptyset$, and ${\cal C}_{n-3}^{i-1}\neq \emptyset$.
\item[(ii)]
 ($\Rightarrow$) Since ${\cal C}_{n-2}^{i-1}={\cal C}_{n-3}^{i-1}=\emptyset$,
  by Lemma~\ref{lemma1}$(ii)$, $i-1>n-2$ or $i-1<\lceil\frac{n-3}{3}\rceil$. If $i-1<\lceil\frac{n-3}{3}\rceil$, then $i-1< \lceil\frac{n-1}{3}\rceil$, and hence ${\cal C}_{n-1}^{i-1}=\emptyset$, a contradiction. So we must have $i>n-1$.
 Also since ${\cal C}_{n-1}^{i-1}\neq \emptyset$, we have $i-1\leq n-1$. Therefore we have $i=n$.\\
($\Leftarrow$) If $i=n$, then by Lemma~\ref{lemma1}$(ii)$,
we have ${\cal C}_{n-2}^{i-1}={\cal C}_{n-3}^{i-1}=\emptyset$ and ${\cal C}_{n-1}^{i-1}\neq \emptyset$.
\item[(iii)]
 ($\Rightarrow$) Since ${\cal C}_{n-1}^{i-1}=\emptyset$, by Lemma~\ref{lemma1}$(ii)$, $i-1>n-1$ or $i-1<\lceil\frac{n-1}{3}\rceil$. If $i-1>n-1$, then $i-1>n-2$ and by lemma~\ref{lemma1}$(ii)$, ${\cal C}_{n-2}^{i-1}= {\cal C}_{n-3}^{i-1}= \emptyset$, a contradiction. So we must have $i<\lceil\frac{n-1}{3}\rceil+1$. But we also have $i-1\geq \lceil\frac{n-2}{3}\rceil$ because ${\cal C}_{n-2}^{i-1}\neq \emptyset$. Hence, we have
$\lceil\frac{n-2}{3}\rceil +1\leq i < \lceil\frac{n-1}{3}\rceil+1$.
Therefore $n=3k+2$ and $i=k+1=\lceil\frac{3k+2}{3}\rceil$ for some $k\in N$.
\\($\Leftarrow$) If $n=3k+2$ and $i=\lceil\frac{3k+2}{3}\rceil$ for some $k\in N$, then by Lemma~\ref{lemma1}$(ii)$, ${\cal C}_{n-1}^{i-1}={\cal C}_{3k+1}^{k}=\emptyset$, ${\cal C}_{n-2}^{i-1}\neq \emptyset$ and ${\cal C}_{n-3}^{i-1}\neq \emptyset$.
\item[(iv)]
 $(\Rightarrow$) Since ${\cal C}_{n-3}^{i-1}= \emptyset$, by Lemma~\ref{lemma1}$(ii)$, we have $i-1>n-3$ or $i-1<\lceil\frac{n-3}{3}\rceil$. Since  ${\cal C}_{n-2}^{i-1}\neq \emptyset$, by Lemma~\ref{lemma1}$(ii)$, we have $\lceil\frac{n-2}{3}\rceil+1\leq i \leq n-1$. Therefore $i-1<\lceil\frac{n-3}{3}\rceil$ is not possible. Hence we must have $i-1>n-3$. Thus $i=n-1$ or $n$. But $i\neq n$ because ${\cal C}_{n-2}^{i-1}\neq\emptyset$. So we have $i=n-1$.
 \\($\Leftarrow$) If $i=n-1$, then by Lemma~\ref{lemma1}$(ii)$, ${\cal C}_{n-1}^{i-1}\neq\emptyset$, ${\cal C}_{n-2}^{i-1}\neq \emptyset$ and ${\cal C}_{n-3}^{i-1}=\emptyset$.
\item[(v)]
 ($\Rightarrow$) Since  ${\cal C}_{n-1}^{i-1}\neq\emptyset,{\cal C}_{n-2}^{i-1}\neq\emptyset$
 and $ {\cal C}_{n-3}^{i-1}\neq \emptyset$, then by applying Lemma~\ref{lemma1}$(ii)$, we have
$\lceil\frac{n-1}{3}\rceil\leq i-1 \leq n-1 , \lceil\frac{n-2}{3}\rceil\leq i-1 \leq n-2$,
 and $ \lceil\frac{n-3}{3}\rceil\leq i-1 \leq n-3.$
So $\lceil\frac{n-1}{3}\rceil\leq i-1 \leq n-3$ and hence $\lceil\frac{n-1}{3}\rceil+1\leq i \leq n-2$.
\\($\Leftarrow$) If $\lceil\frac{n-1}{3}\rceil+1\leq i \leq n-2$, then by Lemma~\ref{lemma1}$(ii)$, we have the result.\quad\qed
\end{enumerate}

\nt The following theorem construct the families of dominating sets of $C_n$.

\begin{theorem}\label{theorem1}  For every $n\geq 4$ and $i\geq \lceil\frac{n}{3}\rceil$,
\begin{enumerate}
\item[(i)]
 If  $\mathcal{C}_{n-1}^{i-1}=\mathcal{C}_{n-2}^{i-1}= \emptyset$ and $\mathcal{C}_{n-3}^{i-1}\neq \emptyset$, then

 $\mathcal{C}_n^i=\mathcal{C}_n^{\frac{n}{3}}=\Big\{\{1,4,\cdots,n-2\},\{2,5,\cdots,n-1\},\{3,6,\cdots,n\}\Big\}$,
\item[(ii)]
If  $\mathcal{C}_{n-2}^{i-1}=\mathcal{C}_{n-3}^{i-1}=\emptyset$ and $\mathcal{C}_{n-1}^{i-1}\neq \emptyset$,
then $\mathcal{C}_n^i=\mathcal{C}_n^n= \Big\{[n]\Big\}$,
\item[(iii)]
 If $ \mathcal{C}_{n-1}^{i-1}=\emptyset$, $\mathcal{C}_{n-2}^{i-1}\neq\emptyset$ and $\mathcal{C}_{n-3}^{i-1}\neq \emptyset$, then
 \\
 $\mathcal{C}_{n}^{i}= \Big\{\{1,4,\cdots,n-4,n-1\},\{2,5,\cdots,n-3,n\},\{3,6,\cdots,n-2,n\}\Big\}\cup
\\ \Big \{ X\cup \left\{
\begin{array}{lll}
\{n-2\},&\mbox{if}& 1\in X \\[10pt]
\{n-1\},&\mbox{if}& 1\not\in X, 2\in X \\[10pt]
\{n\} ,&\mbox{}& otherwise \\[10pt]
\end{array}
|X\in \mathcal{C}_{n-3}^{i-1}\Big \}\right. $

\item[(iv)]
 If $\mathcal{C}_{n-3}^{i-1}=\emptyset, \mathcal{C}_{n-2}^{i-1}\neq\emptyset$ and $\mathcal{C}_{n-1}^{i-1}\neq \emptyset$, then
   $\mathcal{C}_{n}^{i}=\mathcal{C}_{n}^{n-1}=\Big\{[n]-\{x\}|x\in [n]\Big\}$,
 \item[(v)]
  If $\mathcal{C}_{n-1}^{i-1}\neq\emptyset, \mathcal{C}_{n-2}^{i-1}\neq\emptyset$ and $\mathcal{C}_{n-3}^{i-1}\not=\emptyset$, then
 $\mathcal{C}_{n}^{i}=\Big\{\,\{n\}\cup X\;|\;X\in\mathcal{C}_{n-1}^{i-1}\Big\} \cup \\ \Big \{ X_1\cup \left\{
\begin{array}{lll}
\{n\},&\mbox{if}& n-2$ or $ n-3 \in X_1$, for $X_1\in\mathcal{C}_{n-2}^{i-1}\setminus\mathcal{C}_{n-1}^{i-1}  \\[10pt]
\{n-1\},&\mbox{if}&  n-2\not\in X_1, n-3\not\in X_1$ or $X_1\in \mathcal{C}_{n-1}^{i-1}\cap \mathcal{C}_{n-2}^{i-1}\\[10pt]
\end{array}
\Big \}\right.
\cup \\ \Big \{ X_2\cup \left\{
\begin{array}{lll}
\{n-2\},&\mbox{if}& 1\in X_2, $ for $X_2\in \mathcal{C}_{n-3}^{i-1}$ or $X_2\in \mathcal{C}_{n-3}^{i-1}\cap \mathcal{C}_{n-2}^{i-1} \\[10pt]
\{n-1\},&\mbox{if}&  n-3\in X_2$ or $n-4\in X_2$, for $X_2\in \mathcal{C}_{n-3}^{i-1}\setminus \mathcal{C}_{n-2}^{i-1}  \\[10pt]
\end{array}
 \Big \}\right.$.
 \end{enumerate}

 \end{theorem}
\vspace{.1cm}
{\bf Proof.}
\begin{enumerate}
\item[(i)]
 $\mathcal{C}_{n-1}^{i-1}=\mathcal{C}_{n-2}^{i-1}= \emptyset$ and $\mathcal{C}_{n-3}^{i-1}\neq \emptyset$. By Lemma~\ref{lemma4}$(i)$, $n=3k,i=k$ for some $k\in N$. Therefore $\mathcal{C}_n^i=\mathcal{C}_n^{\frac{n}{3}}=\Big\{\{1,4,7,\cdots,n-2\},\{2,5,8,\cdots,n-1\},\{3,6,9,\cdots,n\}\Big\}$.
\item[(ii)]
 $\mathcal{C}_{n-2}^{i-1}=\mathcal{C}_{n-3}^{i-1}=\emptyset$ and $\mathcal{C}_{n-1}^{i-1}\neq \emptyset$. By Lemma~\ref{lemma4}$(ii)$, $i=n$. Therefore $\mathcal{C}_n^i=\mathcal{C}_{n}^{n}=\Big\{[n]\Big\}$.
\item[(iii)]
  $\mathcal{C}_{n-1}^{i-1}=\emptyset, \mathcal{C}_{n-2}^{i-1}\neq\emptyset$, and $ \mathcal{C}_{n-3}^{i-1}\neq \emptyset$. By Lemma~\ref{lemma4}$(iii)$, $n=3k+2, i=k+1$ for some $k\in N$.
   We denote the families $ \Big\{\{1,4,\cdots,3k-2,3k+1\},\{2,5,\cdots,3k-1,3k+2\},\{3,6,\cdots,3k,3k+2\}\Big\}$ and $\Big\{ X\cup \left\{
\begin{array}{lll}
\{3k\},&\mbox{if}& 1\in X \\[10pt]
\{3k+1\},&\mbox{if}& 1\not\in X, 2\in X \\[10pt]
\{3k+2\} ,&\mbox{}& otherwise \\[10pt]
\end{array}
| X\in \mathcal{C}_{3k-1}^{k}\Big\}\right.,
$ by $Y_1$ and $Y_2$, respectively.
We shall prove that $\mathcal{C}_{3k+2}^{k+1}= Y_1 \cup Y_2$. Since  $\mathcal{C}_{3k}^{k}=\Big\{\{1,4,7,\cdots,3k-2\},\{2,5,8,\cdots,3k-1\},\{3,6,9,\cdots,3k\}\Big\}$,
then $ Y_1\subseteq \mathcal{C}_{3k+2}^{k+1}$.
 Also it is obvious that  $Y_2\subseteq \mathcal{C}_{3k+2}^{k+1}$.
 Therefore $Y_1\cup Y_2\subseteq \mathcal{C}_{3k+2}^{k+1}$.
\\
Now let $Y\in \mathcal{C}_{3k+2}^{k+1}$, then by Lemma~\ref{lemma1}$(iii)$,
at least one of the vertices labeled $3k+2,3k+1$ or $3k$ is in $Y$.
 Suppose that $3k+2\in Y$, then by Lemma~\ref{lemma1}$(iii)$, at least one of the vertices labeled $1,2$ or $3$ and $3k+1,3k$ or $3k-1$ are in Y.
If $3k+1$ and at least one of $\{1,2,3\}$, and also $3k$ and at least one of $\{1,2\}$ are in $Y$,
then $Y-\{3k+2\}\in \mathcal{C}_{3k+1}^{k}$, a contradiction.
If $\{3,3k\}$ or $\{2,3k-1\}$ is a subset of $Y$, then $Y=X\cup\{3k+2\}$ for some $X\in \mathcal{C}_{3k}^{k}$.
Hence $Y\in Y_1$. If $\{1,3k-1\}$ is a subset of $Y$, then $Y-\{3k+2\}\in \mathcal{C}_{3k+1}^{k}$, a contradiction.
If $\{3,3k-1\}$ is a subset of $Y$ and $\{3k,3k+1\}$ is not a subset of $Y$, then $Y-\{3k+2\}\in  \mathcal{C}_{3k-1}^{k}$. Hence $Y\in Y_2$.
If $3k+1$ or $3k$ is in $Y$, we also have the result by the similar argument as above.
\item[(iv)]
  By Lemma~\ref{lemma4}$(iv)$, $i=n-1$.
Therefore $\mathcal{C}_n^i=\mathcal{C}_{n}^{n-1}=\Big\{[n]-\{x\}|x\in [n]\Big\}$.
\item[(v)]
 $\mathcal{C}_{n-1}^{i-1}\neq\emptyset, \mathcal{C}_{n-2}^{i-1}\neq\emptyset$ and $\mathcal{C}_{n-3}^{i-1}\neq\emptyset$.
 First, suppose that $X\in\mathcal{C}_{n-1}^{i-1}$, then $X\cup\{n\}\in\mathcal{C}_{n}^{i}$. So $Y_1=\Big\{\,\{n\}\cup X\;|\;X\in \mathcal{C}_{n-1}^{i-1}\Big\}\subseteq \mathcal{C}_{n}^{i}$ .
 Now suppose that $\mathcal{C}_{n-2}^{i-1}\neq \emptyset$. Let $X_1\in\mathcal{C}_{n-2}^{i-1}$. We denote
$ \Big \{ X_1\cup \left\{
\begin{array}{lll}
\{n\},&\mbox{if}& n-2$ or $ n-3 \in X_1$, for $X_1\in \mathcal{C}_{n-2}^{i-1}\setminus\mathcal{C}_{n-1}^{i-1}  \\[10pt]
\{n-1\},&\mbox{if}&  n-2\not\in X_1, n-3\not\in X_1$ or $X_1\in \mathcal{C}_{n-1}^{i-1}\cap \mathcal{C}_{n-2}^{i-1}\\[10pt]
\end{array}
\Big \}\right.$
   simply by $Y_2$.
  By Lemma~\ref{lemma1}$(iii)$, at least one of the vertices labeled $n-3,n-2$ or $1$ is in $X_1$.
  If $n-2$ or $n-3$ is in $X_1$, then $X_1\cup\{n\}\in \mathcal{C}_{n}^i$, otherwise $X_1\cup \{n-1\}\in \mathcal{C}_{n}^{i}$. Hence $Y_2\subseteq \mathcal{C}_{n}^{i}$.
 Here we shall consider $\mathcal{C}_{n-3}^{i-1}\neq \emptyset$. Let $X_2\in \mathcal{C}_{n-3}^{i-1}$. We denote
 $\Big \{ X_2\cup \left\{
\begin{array}{lll}
\{n-2\},&\mbox{if}& 1\in X_2, $ for $X_2\in \mathcal{C}_{n-3}^{i-1}$ or $X_2\in \mathcal{C}_{n-3}^{i-1}\cap \mathcal{C}_{n-2}^{i-1} \\[10pt]
\{n-1\},&\mbox{if}&  n-3\in X_2$ or $n-4\in X_2$, for $X_2\in \mathcal{C}_{n-3}^{i-1}\setminus \mathcal{C}_{n-2}^{i-1}  \\[10pt]
\end{array}
 \Big \}\right.$,
  simply by $Y_3$. If $n-3$ or $n-4$ is in $X$, then $X\cup \{n-1\}\in \mathcal{C}_{n}^{i}$, otherwise $X_2\cup\{n-2\}\in \mathcal{C}_{n}^i$. Hence $Y_3\subseteq Y$. Therefore we've proved that $Y_1\cup Y_2\cup Y_3\subseteq \mathcal{C}_{n}^{i}$.
\\
Now suppose that $Y\in \mathcal{C}_{n}^{i}$, so by Lemma~\ref{lemma1}$(iii)$, $Y$ contain at least one of the vertices labeled $n,n-1$ or $n-2$. If $n\in Y$, so again by Lemma~\ref{lemma1}$(iii)$ at least one of the vertices labeled $n-1,n-2$ or $n-3$  and $1,2$ or $3$ are in $Y$. If $n-2\in Y$ or $n-3\in Y$, then $Y=X\cup\{n\}$ for some $X\in \mathcal{C}_{n-2}^{i-1}$. Hence $Y\in Y_2$.
Otherwise $Y=X\cup \{n-1\}$ for some $X\in \mathcal{C}_{n-2}^{i-1}$.
 Hence $Y\in Y_2$. If $n-1$ or $n-2$ is in $Y$, we also have the result by the similar argument as above.\quad\qed
\end{enumerate}

\vspace{.15cm}

\nt  By Theorem~\ref{theorem1} we have the following theorem for $|\mathcal{C}_{n}^{i}|$.

\vspace{.2cm}

\begin{theorem} \label{theorem2}
If $\mathcal{C}_{n}^{i}$ is the family of dominating set of $C_{n}$ with cardinality $i$, then
\[
|\mathcal{C}_{n}^{i}|=|\mathcal{C}_{n-1}^{i-1}|+|\mathcal{C}_{n-2}^{i-1}|+|\mathcal{C}_{n-3}^{i-1}|.
\]
\end{theorem}
\nt{\bf Proof.} We consider the five cases in Theorem~\ref{theorem1}. We rewrite Theorem~\ref{theorem1} in the following form:
\begin{enumerate}
\item[(i)] If  ${\cal C}_{n-1}^{i-1}={\cal C}_{n-2}^{i-1}= \emptyset$ and ${\cal C}_{n-3}^{i-1}\neq \emptyset$, then
 ${\cal C}_n^i=\\
 \Big\{\{n-2\}\cup X_{1},\{n-1\}\cup X_{2},\{n\}\cup X_{3}|1\in X_{1}, 2\in X_{2}, 3\in X_{3},$ $ X_{1},X_{2},X_{3}\in{\cal C}_{n-3}^{i-1}\Big\}$,
\item[(ii)]
 If  ${\cal C}_{n-2}^{i-1}={\cal C}_{n-3}^{i-1}=\emptyset$ and ${\cal C}_{n-1}^{i-1}\neq \emptyset$, then
 ${\cal C}_n^i=\Big\{\,\{n\}\cup X\;|\;X\in{\cal C}_{n-1}^{i-1}\,\Big\}$,
 \item[(iii)]
  If $ {\cal C}_{n-1}^{i-1}=\emptyset$, ${\cal C}_{n-2}^{i-1}\neq\emptyset$ and ${\cal C}_{n-3}^{i-1}\neq \emptyset$, then
  \\
  ${\cal C}_{n}^{i}=  \Big\{\,\{n\}\cup X_1,\{n-1\}\cup X_2|X_1,X_2\in{\cal C}_{n-2}^{i-1}, 1\in X_2\Big\}\cup
  \\\Big ( X\cup \left\{
\begin{array}{lll}
\{n-2\},&\mbox{if}& 1\in X \\[10pt]
\{n-1\},&\mbox{if}& 1\not\in X, 2\in X \\[10pt]
\{n\} ,&\mbox{}& otherwise \\[10pt]
\end{array}
\Big )\right. $, where $X\in {\cal C}_{n-3}^{i-1}$.
\item[(iv)]
 If ${\cal C}_{n-3}^{i-1}=\emptyset$ and ${\cal C}_{n-2}^{i-1}\neq\emptyset,{\cal C}_{n-1}^{i-1}\neq \emptyset$, then
 \\
 ${\cal C}_{n}^{i}=\Big\{\,\{n\}\cup X_1,\{n-1\}\cup X_2\;|\;X_1\in{\cal C}_{n-1}^{i-1}, X_2\in{\cal C}_{n-2}^{i-1}\, \Big\}$.
\item[(v)]
 If ${\cal C}_{n-1}^{i-1}\neq\emptyset, {\cal C}_{n-2}^{i-1}\neq\emptyset$ and ${\cal C}_{n-3}^{i-1}\not=\emptyset$, then
 \\
 ${\cal C}_{n}^{i}=\Big\{\,\{n\}\cup X\;|\;X\in{\cal C}_{n-1}^{i-1}\Big\} \cup \\ \Big \{ X_1\cup \left\{
\begin{array}{lll}
\{n\},&\mbox{if}& n-2$ or $ n-3 \in X_1$, for $X_1\in{\cal C}_{n-2}^{i-1}\setminus{\cal C}_{n-1}^{i-1}  \\[10pt]
\{n-1\},&\mbox{if}&  n-2\not\in X_1, n-3\not\in X_1$ or $X_1\in {\cal C}_{n-1}^{i-1}\cap {\cal C}_{n-2}^{i-1}\\[10pt]
\end{array}
\Big \}\right.
\cup \\ \Big \{ X_2\cup \left\{
\begin{array}{lll}
\{n-2\},&\mbox{if}& 1\in X_2, $ for $X_2\in {\cal C}_{n-3}^{i-1}$ or $X_2\in {\cal C}_{n-3}^{i-1}\cap {\cal C}_{n-2}^{i-1} \\[10pt]
\{n-1\},&\mbox{if}&  n-3\in X_2$ or $n-4\in X_2$, for $X_2\in {\cal C}_{n-3}^{i-1}\setminus{\cal C}_{n-2}^{i-1}  \\[10pt]
\end{array}
 \Big \}\right.$.
  where $X_1\in {\cal C}_{n-2}^{i-1}\setminus {\cal C}_{n-1}^{i-1}$ and $X_2\in {\cal C}_{n-3}^{i-1}\setminus {\cal C}_{n-2}^{i-1}\cap {\cal C}_{n-1}^{i-1}$.
\\
By above construction, in every cases, we have $|{\cal C}_{n}^{i}|=|{\cal C}_{n-1}^{i-1}|+|{\cal C}_{n-2}^{i-1}|+|{\cal C}_{n-3}^{i-1}|.$\quad\qed
\end{enumerate}

\nt Since $|{\cal C}_n^i|$ satisfy the recursive formula with two variable, finding a formula for $|{\cal C}_n^i|$ is not easy.
 In the following theorem  we use the generating function technique to find $\Big|{\cal C}_n^i\Big|$.

\begin{theorem}
For every natural $n\geq 4$ and $\lceil\frac{n}{3}\rceil\leq i \leq n$, $|{\cal C}_n^i|$
is the coefficient of $u^nv^i$ in the expansion of the function
\[
 f(u,v)=\frac{u^4v^2(6+4v+v^2+3u+4uv+uv^2+u^2+3u^2v+u^2v^2)}{1-uv-u^2v-u^3v}.
 \]
\end{theorem}
\nt{\bf Proof.} Set $f(u,v)=\sum_{n=4}^\infty\sum_{i=2}^\infty |{\cal C}_n^i|u^nv^i$.
By recursive formula for $|{\cal C}_n^i|$ in  Theorem~\ref{theorem2} we can write $f(u,v)$ in the following form
\[ f(u,v)=\sum_{n=4}^\infty \sum_{i=2}^\infty (|{\cal C}_{n-1}^{i-1}|+|{\cal C}_{n-2}^{i-1}|+|{\cal C}_{n-3}^{i-1}|)u^nv^i=
\]
\[
uv\sum_{n=4}^\infty \sum_{i=2}^\infty |{\cal C}_{n-1}^{i-1}|u^{n-1}v^{i-1}+u^2v\sum_{n=4}^\infty \sum_{i=2}^\infty |{\cal C}_{n-2}^{i-1}|u^{n-2}v^{i-1}+\]
\[u^3v\sum_{n=4}^\infty \sum_{i=2}^\infty |{\cal C}_{n-3}^{i-1}|u^{n-3}v^{i-1}=uv(|{\cal C}_3^1|u^3v+|{\cal C}_3^2|u^3v^2+|{\cal C}_3^3|u^3v^3)+ uvf(u,v)+\]
\[+u^2v(|{\cal C}_2^1|u^2v+|{\cal C}_2^2|u^2v^2+|{\cal C}_3^1|u^3v+|{\cal C}_3^2|u^3v^2+|{\cal C}_3^3|u^3v^3)
+u^2vf(u,v)+\]
\[u^3v( |{\cal C}_1^1|uv+|{\cal C}_2^1|u^2v+|{\cal C}_2^2|u^2v^2+|{\cal C}_3^1|u^3v+|{\cal C}_3^2|u^3v^2+|{\cal C}_3^3|u^3v^3)
+u^3vf(u,v)
\]
By substituting the values from Table 1, we have
\[
f(u,v)(1-uv-u^2v-u^3v)=u^4v^2(6+4v+v^2+3u+4uv+uv^2+u^2+3u^2v+u^2v^2)
\]
Therefore we have the result.\quad\qed

\section{Domination polynomial of a cycle}
In this section we introduce and investigate the domination polynomial of cycles.
\begin{define}\label{Definition1}
 Let ${\cal C}_n^i$ be the family of dominating sets of a cycle $C_n$ with cardinality $i$ and let
$d(C_n,i)=|{\cal C}_n^i|$. Then the domination polynomial $D(C_n,x)$ of $C_n$ is defined as
\begin{center}
$D(C_n,x)=\sum_{i=\lceil\frac{n}{3}\rceil}^{n} d(C_n,i) x^{i}$.
\end{center}

\end{define}
\nt By Definition~\ref{Definition1} and Theorem~\ref{theorem2}, we have the following theorem.
\begin{theorem} \label{}
For every $n\geq 4$,
\[
D(C_{n},x)=x\Big[D(C_{n-1},x)+D(C_{n-2},x)+D(C_{n-3},x)\Big],
\]
 with the initial values $D(C_{1},x)=x, D(C_{2},x)=x^{2}+2x, D(C_{3},x)=x^{3}+3x^{2}+3x$.\quad\qed
\end{theorem}

\nt Using  Theorem~\ref{theorem2}, we obtain the coefficients of $D(C_{n},x)$ for $1\leq n\leq 16$ in Table 1.
Let $d(C_n,j)=|{\cal C}_n^j|$.
 There are interesting relationships between the numbers $d(C_n,j) (\frac{n}{3}\leq j\leq n)$
in the table 1.

\[
\begin{footnotesize}
\small{
\begin{tabular}{r|lcrrrcccccccccccc}
$j$&$1$&$2$&$3$&$4$&$5$&$6$&$7$&$8$&$9$&$10$&$11$&$12$&$13$&$14$&$15$&$16$\\[0.3ex]
\hline
$n$&&&&&&&&&&&&&&&&\\
$1$&1&&&&&&&&&&&&&&&\\
$2$&2&1&&&&&&&&&&&&&&\\
$3$&3&3&1&&&&&&&&&&&&&\\
$4$&0&6&4&1&&&&&&&&&&&&\\
$5$&0&5&10&5&1&&&&&&&&&&&\\
$6$&0&3&14&15&6&1&&&&&&&&&&\\
$7$&0&0&14&28&21&7&1&&&&&&&&&\\
$8$&0&0&8&38&48&28&8&1&&&&&&&&\\
$9$&0&0&3&36&81&75&36&9&1&&&&&&&\\
$10$&0&0&0&25&102&150&110&45&10&1&&&&&&\\
$11$&0&0&0&11&99&231&253&154&55&11&1&&&&&\\
$12$&0&0&0&3&72&282&456&399&208&66&12&1&&&&\\
$13$&0&0&0&0&39&273&663&819&598&273&78&13&1&&&\\
$14$&0&0&0&0&14&210&786&1372&1372&861&350&91&14&1&&\\
$15$&0&0&0&0&3&125&765&1905&2590&2178&1200&440&105&15&1&\\
$16$&0&0&0&0&0&56&608&2214&4096&4560&3312&1628&544&120&16&1&\\
\end{tabular}}
\end{footnotesize}
\]
\begin{center}
\nt{Table 1.} $d(C_{n},j)$ The number of dominating sets of $C_n$ with cardinality $j$.
\end{center}

\nt  In the following theorem, we obtain some properties of $d(C_n,j)$:
\begin{theorem} \label{theorem9}
The following properties hold for coefficients of $D(C_n,x)$:
\begin{enumerate}
\item[(i)]
 For every $n\in N$, $d(C_{3n},n)=3$,
\item[(ii)]
 For every $n\geq 4, j\geq \lceil\frac{n}{3}\rceil$, $d(C_{n},j)=d(C_{n-1},j-1)+d(C_{n-2},j-1)+d(C_{n-3},j-1)$,
\item[(iii)]
 For every $n\in N$, $d(C_{3n+2},n+1)=3n+2$,
\item[(iv)]
 For every $n\in N$, $d(C_{3n+1},n+1)=\frac{n(3n+7)+2}{2}$,
\item[(v)]
 For every $n\geq 3$, $d(C_n,n)=1$,
\item[(vi)]
 For every $n\geq 3$, $d(C_n,n-1)=n$,
\item[(vii)]
 For every $n\geq 3$, $d(C_n,n-2)=\frac{(n-1)n}{2}$,
\item[(viii)]
 For every $n\geq 4$, $d(C_n,n-3)=\frac{(n-4)(n)(n+1)}{6}$,
\item[(ix)]
 for every $j\geq 4$, $ \sum_{i=j}^{3j}d(C_{i},j)=3\sum_{i=j-1}^{3j-3}d(C_{i},j-1),$
\item[(x)]
 for every $n\geq 3$, $ 1=d(C_{n},n)< d(C_{n+1},n)< d(C_{n+2},n)<\cdots<d(C_{2n-1},n)<d(C_{2n},n)>d(C_{2n+1},n)>\cdots>d(C_{3n-1},n)>d(C_{3n},n)=3.$
 \item[(xi)]
 If $S_n=\sum_{j=\lceil\frac{n}{3}\rceil}^n d(C_n,j)$, then for every $n\geq 4$, $S_n=S_{n-1}+S_{n-2}+S_{n-3}$
with initial values $S_1=1, S_2=3$ and $S_3=7$.
\end{enumerate}
\end{theorem}
\nt{\bf Proof.}
\begin{enumerate}
\item[(i)]
 Since  ${\cal C}_{n}^{3n}=\Big\{\{1,4,7,...,3n-2\},\{2,5,8,...,3n-1\},\{3,6,9,...,3n\}\Big\}$, so \\$d(C_{3n},n)=3$.
\item[(ii)]
 It follows from Theorem~\ref{theorem2}.
\item[(iii)]
 By induction on $n$. The result is true for $n=1$, because \\${\cal C}_2^5=\Big\{\{1,3\},\{1,4\},\{2,4\},\{2,5\},\{3,5\}\Big\}$.
 Now suppose that the result is true for all natural numbers less than $n$, and we prove it for $n$. By $(i),(ii)$ and induction hypothesis, we have
\begin{eqnarray*}
d(C_{3n+2},n+1)&=&d(C_{3n+1},n)+d(C_{3n},n)+d(C_{3n-1},n)\\
&=&3n+2.
\end{eqnarray*}
\item[(iv)]
 By induction on $n$. Since ${\cal C}_{2}^4=\Big\{ \{ 1,2\},\{1,3\},\{1,4\},\{2,3\},\{2,4\},\{3,4\}\Big\}$,
 so $d(C_4,2)=6$, the result is true for $n=1$.
  Now suppose that the result is true for all natural numbers less than $n$, and we prove it for $n$.
 By $(i),(ii),(iii)$ and induction hypothesis, we have
\begin{eqnarray*}
d(C_{3n+1},n+1)&=&d(C_{3n},n)+d(C_{3n-1},n)+d(C_{3n-2},n)\\
&=&3+3(n-1)+2+\frac{(n-1)(3(n-1)+7)+2}{2}\\
&=&\frac{n(3n+7)+2}{10}.
\end{eqnarray*}
\item[(v)] Since for any graph with $n$ vertices, $d(G,n)=1$, then we have the result.
\item[(vi)]
 Since for any graph with $n$ vertices, $d(G,n-1)=n$, we have the result.
\item[(vii)]
 By induction on $n$. The result is true for $n=3$, since $d(C_3,1)=3$.
 Suppose that the result is true for all natural number less than $n$, and we prove it for $n$. By parts $(ii),(v),(vi)$ and induction hypothesis have
\begin{eqnarray*}
d(C_n,n-2)&=&d(C_{n-1},n-3)+d(C_{n-2},n-3)+d(C_{n-3},n-3)\\
&=&
\frac{(n-2)(n-1)}{2}+n-2+1\\
&=&\frac{(n-1)n}{2}
\end{eqnarray*}
\item[(viii)]
 By induction on $n$. The result is true for $n=4$, since $d(C_4,1)=0$.
 Suppose that the result is true for all natural number less than $n$ and prove it for $n$. By parts $(ii),(vi),(vii)$ and induction hypothesis we have
\begin{eqnarray*}
d(C_n,n-3)&=&d(C_{n-1},n-4)+d(C_{n-2},n-4)+d(C_{n-3},n-4)\\
&=&\frac{(n-5)(n-1)n}{6}+\frac{(n-2)(n-3)}{2}+n-3\\
&=&\frac{(n-4)n(n+1)}{6}
\end{eqnarray*}
\item[(ix)]
 Proof by induction on $j$. First, suppose that $j=3$. Then $\sum_{i=3}^{9}d(C_{i},3)=54=3\sum_{i=2}^{6}d(C_{i},2)$.
 Now suppose that the result is true for every $j<k$, and we prove for $j=k$:
\begin{eqnarray*}
\sum_{i=k}^{3k}d(C_{i},k)&=& \sum_{i=k}^{3k} d(C_{i-1},k-1)+\sum_{i=k}^{3k}d(C_{i-2},k-1)+\sum_{i=k}^{3k}d(C_{i-3},k-1)\\
&=&
3\sum_{i=k-1}^{3(k-1)}d(C_{i-1},k-2)+3\sum_{i=k-1}^{3(k-1)}d(C_{i-2},k-2)\\
&+&3\sum_{i=k-1}^{3(k-1)}d(C_{i-3},k-2)=3\sum_{i=k-1}^{3k-3}d(C_{i},k-1).
\end{eqnarray*}
\item[(x)]
 We shall prove that for every $n$, $d(C_{i},n)<d(C_{i+1},n)$ for $n\leq i \leq 2n-1$,
 and $d(C_{i},n)>d(C_{i+1},n)$ for $2n\leq i \leq3n-1$. We prove the first inequality by induction on $n$.
 The result hold for $n=3$.
 Suppose that result is true for all $n\leq k$. Now we prove it for $n=k+1$, that is $d(C_{i},k+1)<d(C_{i+1},k+1)$ for $k+1\leq i \leq 2k+1$. By Theorem~\ref{theorem2}
 and induction hypothesis we  have
\begin{eqnarray*}
d(C_{i},k+1)&=&d(C_{i-1},k)+d(C_{i-2},k)+d(C_{i-3},k)\\
&<&d(C_{i},k)+d(C_{i-1},k)+d(C_{i-2},k)\\
&=&d(C_{i+1},k+1)
\end{eqnarray*}
Similarly, we have the other inequality.
\item[(xi)]
By Theorem~\ref{theorem2}, we have
\begin{eqnarray*}
S_n&=&\sum_{j=\lceil\frac{n}{3}\rceil}^n d(C_n,j)=\sum_{j=\lceil\frac{n}{3}\rceil}^n(d(C_{n-1},{j-1})+d(C_{n-2},{j-1})+d(C_{n-3},j-1))\\
&=&
\sum_{j=\lceil\frac{n}{3}\rceil-1}^{n-1}d(C_{n-1},j)+\sum_{j=\lceil\frac{n}{3}\rceil-1}^{n-2}d(C_{n-2},j)
+\sum_{j=\lceil\frac{n}{3}\rceil-1}^{n-3}d(C_{n-3},j-1)\\
&=&
S_{n-1}+S_{n-2}+S_{n-3}.\quad\qed
\end{eqnarray*}

\end{enumerate}

\end{document}